\magnification=\magstep1

\newcount\sec \sec=0
\input Ref.macros
\input math.macros
\input labelfig.tex
\forwardreferencetrue \citationgenerationtrue \initialeqmacro
\sectionnumberstrue
\input epsf


\def\dist{{\rm dist}}


\title{Invariant colorings of random planar maps}

\author{\'Ad\'am Tim\'ar}
\bigskip

\abstract{We show that every locally finite random 
graph embedded in the plane with an isometry-invariant distribution can be 
5-colored in an invariant and 
deterministic way, under some nontriviality assumption and a mild 
assumption on the tail of edge lengths. The assumptions hold for any 
Voronoi map on a point process that has no nontrivial symmetries almost 
surely, hence we improve and generalize previous results on 
6-coloring the Voronoi map on a Poisson point process \ref b.ABGMP/. 
}

\bsection{Introduction}{s.1}

We consider random graphs $G$ embedded in the plane such that the number 
of vertices in any bounded set is a.s. finite, and such that 
the 
distribution of the image of $G$ as a subset of the plane (which we also 
denote by $G$) is {\bf 
invariant} with respect 
to some transitive group $\Gamma$ of isometries of the plane (e.g. all 
isometries, 
or all the 
translations). With a slight abuse of notation, we use $G$ both for the 
graph and for the embedded image of it in the plane (hence thinking about 
it as a 1-complex in the plane). We want to give a coloring of this graph 
with as few colors 
as possible, in such a way that the coloring is an equivariant and 
measurable function of 
$G$ with respect to $\Gamma$. In other words, if $B_R(x)$ is the disc of 
radius $R$ around $x$ in the plane,
we 
would like to construct a mapping $c_G :V(G)\to \{1,\ldots ,k\}$ such 
that $c_G$ assigns different values to adjacent vertices, and $c_G$ 
satisfies: 
\item{(i)} $c_{g (G)} (g (x))= c_G (x)$ for every
$g\in\Gamma$,
\item{(ii)} with probability tending to 1 with $R$, $c_G (x)$ can be 
determined from 
$B_R(x)\cap G$,
as a measurable function of $B_R(x)\cap G$.
 
We want to make $k$ as small as possible.
The condition of ``being an equivariant 
function of $G$" is sometimes simply referred to as ``equivariance". 
Note that the assumption on measurability is important not only in order 
to avoid a trivial answer (since every infinite planar graph is 
4-colorable, as we will discuss later), but also in order to be able to 
talk about probabilistic properties of the coloring (e.q. the 
distribution of what one sees locally around some fixed point of the 
underlying space). 
This natural 
assumption is always needed when one studies 
equivariant functions of point processes (even though in some papers this 
assumption is not made explicitly). In the paper for simplicity when we 
say that a function is equivariant, 
we 
always include that it is also measurable.

We will impose the condition that the set of symmetries of $G$ is 
almost always trivial, where by a {\bf symmetry} we mean a graph 
automorphism 
which is achieved by an isometry in the plane. In particular, the 
condition holds for 
any random graph $G$ where the graph has no nontrivial automorphism, 
or for $G$ such that the only isometry of $V(G)$ (as a point set in the 
plane) is the trivial one. We will also need a condition on $G$, which is 
in brief a condition on the 
existence of relatively few long edges. 
Namely, say that a $G$ invariant planar map has the {\bf regular decay 
property}  
if
for every $r$ there is an $a(r)$, with $a(r)\to 0$ as $r\to \infty$, such
that for any $R\geq r$, the
probability
that the $R$-neighborhood $B_R$ of $o$ in the plane
contains the endpoints
$u,v$ of an edge in $G$ such that $\dist_{\R^2} (u,v)\geq a(r) R/6$
is smaller than $a(r)$.

\procl t.main
Let $G$ be a random graph in the plane, whose distribution is invariant 
with respect to some transitive group of isometries of the plane. Suppose 
that with probability 1, the only symmetry that $G$ has is the 
trivial one and that $G$ satisfies the regular decay property. Then there 
exists a 5-coloring for $G$ which is an 
equivariant 
function of $G$.

In particular, every Voronoi map $G$ on some point process with only the 
trivial symmetry has a 5-coloring which is an equivarint function of the 
point process. 
\endprocl

See \ref r.lattice/ for the case when the symmetry assumption on $G$ does 
not hold, in which case $G$ is a ``quasi" lattice with finitely many 
orbits. Such graphs either trivially have no equivariant coloring by any 
finite 
number of colors, or their ``equivariant measurable" chromatic number is 
7 or less 
(depending on the chromatic number of the factor graph $G/\Gamma$). 

The famous 4 color theorem states that every finite planar map is 
4-colorable. This was first proved by Appel and Haken, then 
with a much smaller, but still significant amount of computer verification 
by 
Robertson, Sanders, Seymour and Thomas. (See \ref b.Th/ for a survey on 
the 4 color theorem, and further references.) This can be extended to 
infinite graphs, 
by standard compactness arguments. However, if one did this right away, 
one is likely to get a 4-coloring that is neither equivariant, nor 
measurable.



The question that we address in \ref t.main/ was asked by Itai Benjamini, 
and an equivariant coloring by 6 colors is given for Voronoi 
tessellation on a Poisson point process by Angel, Benjamini, 
Gurel-Gurevich, Mayerovich, Peled \ref b.ABGMP/. The proof in 
\ref b.ABGMP/ uses 
some explicit computations about the distribution of the number of 
neighbors of a region in this planar map, and the bounds attained are 
used to 
prove that by repeatedly removing every region of $\leq 5$ neighbors from 
the graph, one gets only finite components, after a finite number of 
iterations.
This need not be true for a general $G$ that only satisfies the assumption 
in \ref 
t.main/ (even not for every Voronoi map on a point process, see \ref 
x.uj/), 
hence the 
proof in \ref b.ABGMP/ (whose second part is a ``greedy" coloring, as in 
the usual proof 
of 
the 6 color theorem), does not seem to fully generalize to our setup, even 
with 
6 
colors.

\ref t.main/ will follow from \ref l.cycles/ and \ref l.exhaustion/ right 
away. For some further relaxation of the condition on $G$, see \ref 
r.generalizations/. 

The reason we need that there are no nontrivial symmetries is
that then there is a so called {\bf index function} from $G$ 
(respectively from $G\times
G$) to the reals, which is an {\it injective} equivariant (respectively 
diagonally invariant) function of
$G$. This
enables one to take certain subsets of infinite point sets in an 
equivariant 
way, or make {\it local choices}, e.g. choose a vertex from each finite
class of some equivariant partition,
and still preserve equivariance. See \ref b.T/ for more details.
Hence, when we say ``choose", ``fix" etc. some point of each element of 
some equivariant collection of finite subsets
of $G$, it always means that we 
have some 
previously fixed rule, which makes the choice depend on the precise 
local configuration (using the index function), and makes it remain 
equivariant and a deterministic 
function of the configuration.


An {\bf induced subgraph} of $G$ is a subgraph $H$ such that the set of 
edges of $G$
with both endpoint in $V(H)$ is equal to $H$. We call a subgraph of $G$ 
{\bf
non-selftouching}, if the graph that it induces in $G$ is itself. Two 
subgraphs of $G$ are called {\bf non-touching}, if they are not adjacent: 
there 
is no 
edge with one endpoint in each.
By a {\bf path} we always mean a simple path, i.e. no multiple vertices 
are allowed.

In \ref s.2/ we present the combinatorial trick that reduces the question 
of coloring to finding a certain kind of exhaustion for $G$. The existence 
of such an exhaustion is less sensitive to local changes than colorings. 
\ref s.3/ contains the construction of such an exhaustion, with some 
complications because of the generality of our setup. \ref s.4/ concerns 
some open questions and generalizations. 

\def\int{{\rm int}}

\bsection{5-coloring from induced cycles}{s.2}
Given a cycle $C$ in $G$, define $\int (C)$ to be the subgraph induced 
in $G$ by the set of vertices 
in the bounded component of $G\setminus C$. For the next definition, 
note that for any infinite tree with one end, one can define a parent 
to each vertex $w$, as the first vertex on the path from 
$w$ to 
infinity. If $v$ is the parent of $w$, we will write $w\to v$.

Say that $(T,\lambda)$ is an {\bf even cycle exhaustion} of $G$, {\bf 
with corridors of width} $c>0$, if it is 
an equivariant function of $G$, and 
\item{(1)} $T$ is an infinite 
tree with one-end; 
\item{(2)} $\lambda : V(T)\mapsto 2^G$ is such that for every vertex $v$ 
of 
$T$, 
$\lambda (v)$ is a non-selftouching cycle of {\it even} length of $G$; 
\item{(3)} any $\lambda (v)$ and $\lambda (w)$ have distance $\geq c$ 
whenever 
$v\not = w$;
\item{(4)} $\lambda (w)\subset \int(\lambda (v))$ whenever $w\to v$.
\item{(5)} $G\setminus \cup_{v\in V(T)} \lambda (v)$ has only finite 
components.

If we do not require the $\lambda (v)$'s to have even lengths, then we 
simply call the above structure a {\it cycle exhaustion}.

Informally, an even cycle exhaustion is a collection of non self-touching 
cycles 
of 
even lengths, such that their pairwise distances are at least $c$, every 
point of the plane is surrounded by at least one (and hence infinitely 
many) of these cycles, and the relation ``surrounding" defines a 
natural tree-structure on these cycles. 
Note that by defining the collection of cycles, the tree 
structure is uniquely defined as well.
This is how we prefer to think 
about $(T,\lambda )$. 

\procl l.cycles
Let $(T,\lambda)$ be an even cycle exhaustion of $G$ of corridor width 
$4$. Then there is 
an 
equivariant 5-coloring of $G$.
\endprocl

\def\w_v{p_w (v)}
\def\v_v{p_v (v)}
\def\w_w{p_w (w)}
\def\v_w{p_v (w)}

\proof
First, for every cycle $\lambda (v)$, which is a bipartite cycle, fix one 
of the two classes of the bipartition, and say that its elements are the 
odd elements, while the elements in the other class are the even elements. 
Do it so that the choices are invariant with respect to $\Gamma$.

For every $v\in V(T)$ consider the finite graph $H_v$ induced by 
$\int(\lambda (v))\cup\lambda (v)$ in G. Then in $H_v$, for every $w\to 
v$, contract 
$H_w$ (which is naturally sitting in $H_v$). Call the resulting vertex 
$p_v (w)$. Also contract $\lambda (v)$ to one vertex $p_v (v)$ in $H_v$. 
The 
graph obtained from $H_v$ after these contractions is called $G_v$. The 
vertices of $G_v$ that did not come to existence by contraction, but were 
present in $H_v$, are called {\bf ordinary} vertices.
See Figure 1 for an example (note however, that the example there does 
not satisfy the condition on the corridor width).
We refer to
ordinary vertices 
and their identical copy in $G$ under the same name.
Since $G_v$ is finite planar, we can fix some 4-coloring 
$\gamma_v :V(G_v)\mapsto\{1,2,3,4\}$
of $G_v$, making the choices invariant under $\Gamma$. 
$$$$
\SetLabels
  (.2*1.03) {$\lambda (v)$}\\
 (.2*.71) {$\lambda (w)$}\\
 (.6*.03) {$p_v (v)$}\\
 (.87*.44) {$p_v (w)$}\\
  (.5*-.17) {{\bf Figure 1.} $H_v$ (left) and $G_v$ (right). We simplified 
the picture by making $v$ 
have}\\
  (.5*-.26) {only one child, and the distance between $\lambda (v)$ and 
$\lambda (w)$ be only 2.}\\ 
\endSetLabels
\AffixLabels{\epsfysize=4.8cm
\epsfbox{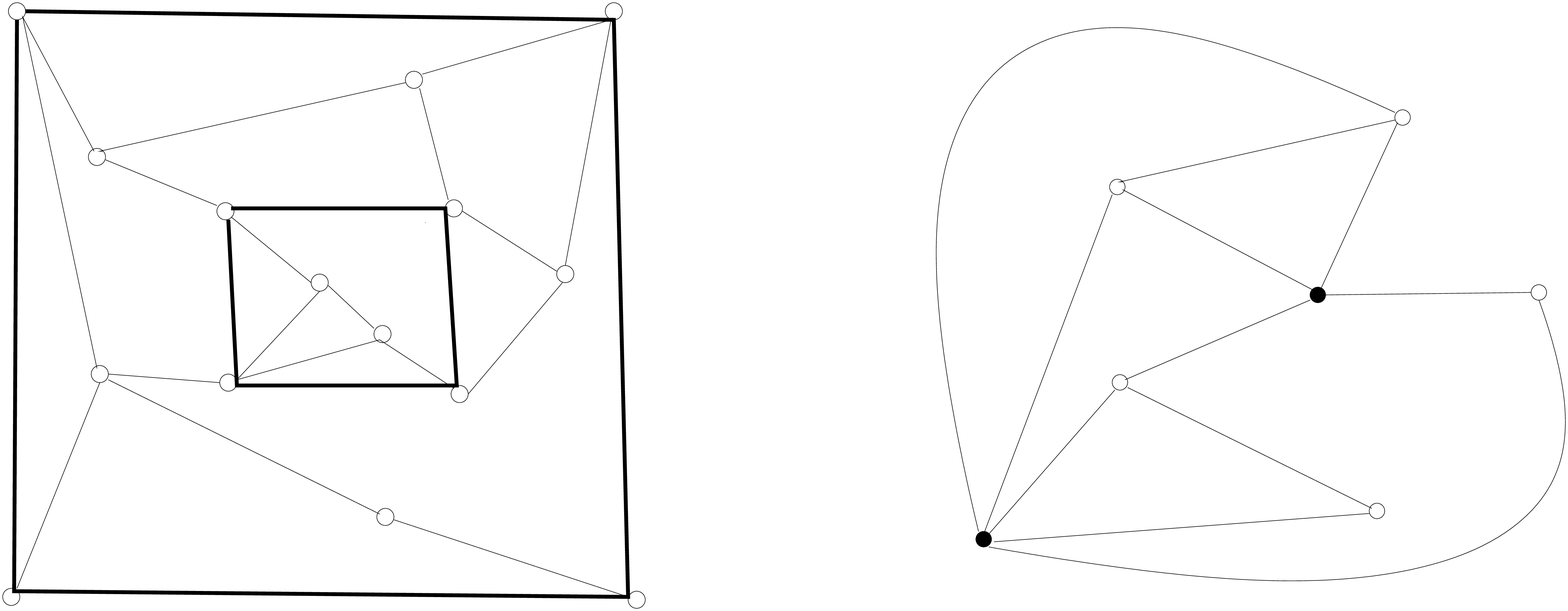}
}
$$\,$$
$$$$

To get a coloring of $G$, do the following. For every $v\in V(t)$, 
and every ordinary
$x\in G_v$ 
assign $x$ the color given to it by 
$\gamma_v$. This coloring of $G\setminus \cup_{v\in V(T)}\lambda 
(v)$ will be called $\gamma$. 

Now, every cycle $\lambda (w)$ was contracted into vertex $\w_w$ in 
$G_w$, and into vertex $\v_w$ in $G_v$, where $w\to v$. If $\gamma_v 
(\v_w )=\gamma_w (\w_w)$, color every even vertex of $\lambda (w)$ with 
color 
$\gamma_v (\v_w )$, and every odd vertex with color 0. 
See Figure 2 for an illustration of this case.

Otherwise, if  $\gamma_v    
(\v_w )\not =\gamma_w (\w_w )$, color every even vertex of $\lambda (w)$ 
with $\gamma_v (\v_w )$, and every odd one with $\gamma_w(\w_w )$.

Call 
the 
resulting assignment of colors to $V(G)$ (which we obtain 
by extending $\gamma$ from $G\setminus \cup\lambda (w)$ to $G$ as just 
described) $\gamma'$. This $\gamma'$ is 
typically 
not a good coloring yet. There may be $G$-neighbors of identical color in 
two possible ways: either an element of $\lambda (w)$ and 
an ordinary vertex in $G_w$ both got color $\gamma_v (\v_w )$, or 
an element of 
$\lambda (w)$ and an ordinary vertex in $G_v$ both got color $\gamma_w 
(\w_w )$.
For all such pairs, recolor the point not in $\gamma (w)$ by assigning 
it color 0. Doing 
this for all pairs of neighbors that had the same color by $\gamma'$, we 
obtain a coloring $\Gamma$ of $V(G)$, which we claim to be a good 
coloring. For this, one only has to check that a vertex $x$ that was 
recolored to 0 in this last step, has no recolored neighbor $y$, and 
no neighbor $z$ that had color 0 by $\gamma'$. 
The existence of a $z$ as above is not possible by the condition that  
the $\gamma (w)$'s are at distance $\geq 5$ from each other, and every
recolored vertex is at distance 1 from some $\gamma (w)$.
The existence of an $y$ as above is not possible because, 
if there existed such a $y$, then one would have $\gamma 
(x)=\gamma (y)$, and thus $\gamma_u (x)=\gamma_u (y)$ with the 
appropriate 
$u$ ($x,y\in G_u$), which 
would contradict that $\gamma_u$ is a good coloring. See Figure 3 for this 
case.
This finishes the proof that $\Gamma$ is a 
5-coloring as desired.\Qed

$$$$
\SetLabels
  (.2*.93) {$\lambda (v)$}\\
 (.2*.71) {$\lambda (w)$}\\
 (.68*.24) {$p_w (w)$}\\
 (.87*.63) {$p_v (w)$}\\
 (.56*.35) {$p_v (v)$}\\
  (.53*-.17) {{\bf Figure 2.} The coloring of $\lambda (v)\cup\int
(\lambda (v))$
(left) coming from $\gamma_v$ of $G_v$ (right, upper), and}\\
  (.5*-.24) {$\gamma_w$
for $G_w$ (right, lower), when $\gamma_v (p_v (v))=\gamma_w (p_w
(v))$. Black stands for color 0.}\\
\endSetLabels
\AffixLabels{\epsfysize=6.2cm
\epsfbox{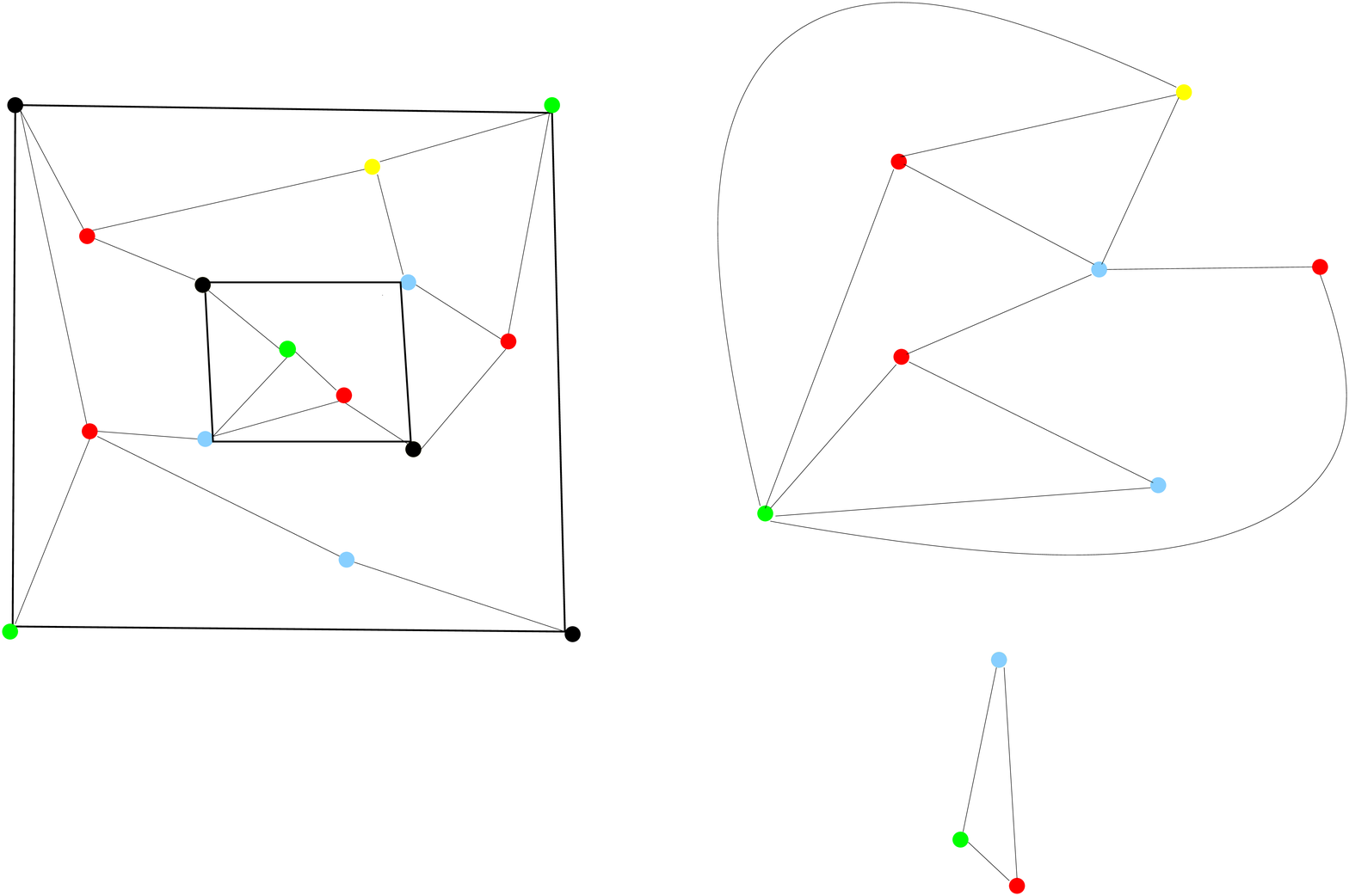}}
$$\,$$
$$\,$$

$$$$
\SetLabels
  (.53*-.17) {{\bf Figure 3.} The coloring of $\int (\lambda (v))$
(left) coming from $\gamma_v$ on $G_v$ (right, upper),}\\
  (.5*-.24) {and $\gamma_w$
for $G_w$ (right, lower), when $\gamma_v (p_v (v))\not=\gamma_w (p_w 
(v))$.}\\
\endSetLabels
\AffixLabels{\epsfysize=6.2cm
\epsfbox{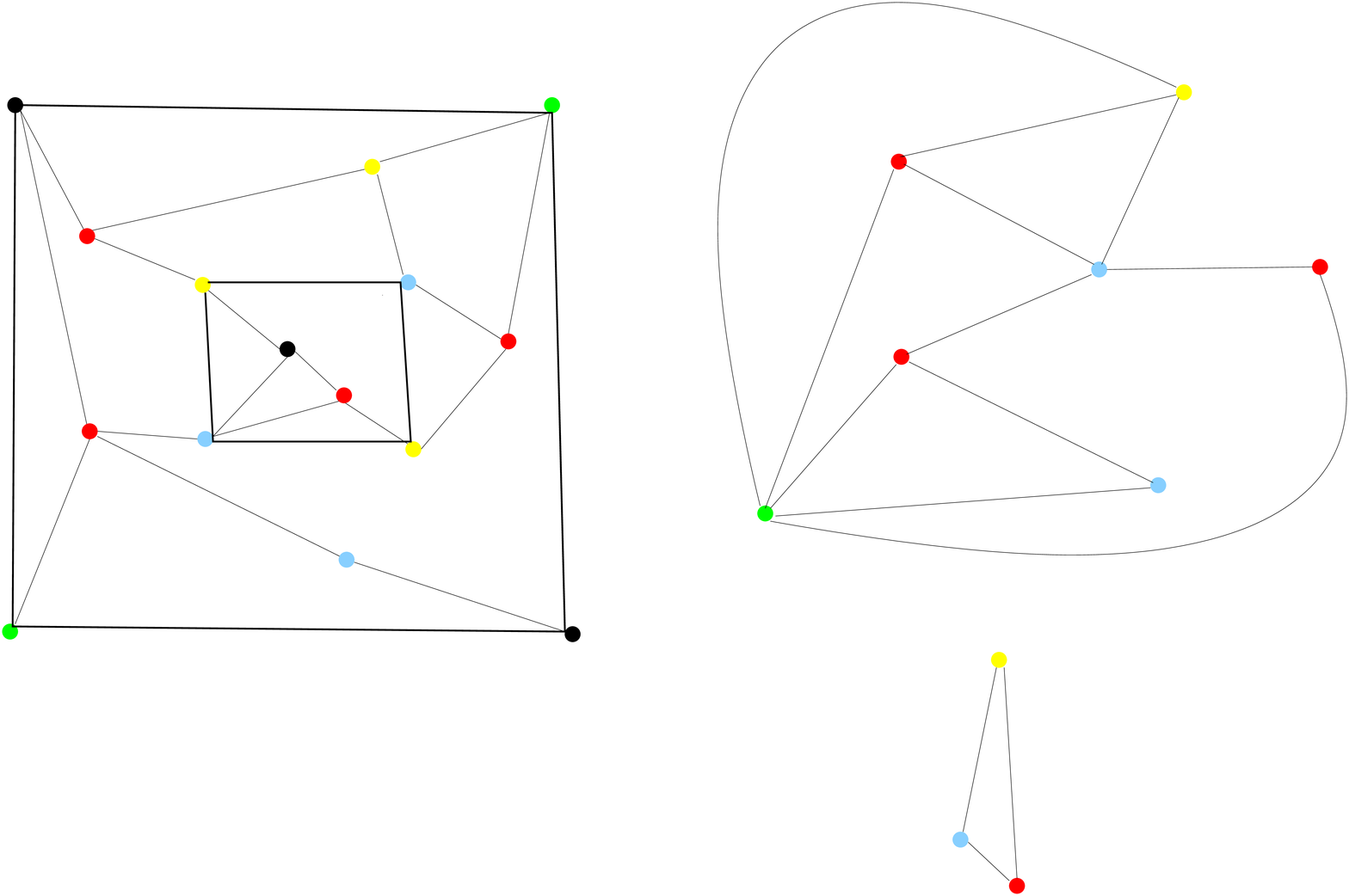}}   
$$$$
\bigskip
\medskip


\bsection{Existence of a cycle exhaustion}{s.3}

\def\pi{{\cal P}_i}
We shall {\it assume that} $G$ {\it is triangulated}. This is not a 
restriction: we can triangulate every face of $G$ in some equivariant 
deterministic way, also respecting the conditions on 
$G$. An equivariant coloring of the new, triangulated graph is
also a coloring for the original one.

Given a subgraph $H$ of $G$, let $\partial H$ be the outer boundary of 
$H$, that is, the set of vertices 
in $G\setminus H$ that are adjacent to $H$. Let $\partial_r H$ be the 
set of vertices in $G\setminus H$ at distance at most $r$ from $H$.

We will need the following graph theoretic observation later:

\procl p.boundary
Let $H$ be some connected subgraph of the infinite, triangulated planar 
graph $G$. Then the set $O(H)$ 
of vertices $x$ in $\partial H$ that are visible from infinity (i.e. there 
is an infinite path from $x$ in $(G\setminus ( \partial H\setminus  
\{x\})$) 
induces a non-selftouching cycle in $G$.
\endprocl

\proof
There is a natural cyclic ordering on $O(H)$ (defined as we ``walk 
along" $O(H)$ in $G\setminus (H\cup\partial H)$, and 
any two vertices 
following each other in this ordering are adjacent, because $G$ is 
triangulated. Thus there is a cycle with vertex set $O(H)$, and we only 
have to prove that $O(H)$ induces no edges other than these. Now, if 
$O(H)$ induced some other edge $\{x,y\}$, 
then the graph induced in $G$ by $O(H)\setminus \{x,y\}$ would have at 
least two 
components, because $x$ and $y$ 
are both visible from $H$ and from infinity. This contradicts the 
fact that any two vertices of $O(H)$ can be joined by a path with 
every inner vertex in $H$, which should be true since $H$ is connected 
and $O(H)$ is in its boundary.
\Qed

Call the set of vertices in $\partial H$ visible from infinity the
{\bf exterior boundary} of $H$.

A much stronger version of the next lemma was proved in \ref b.T2/, for 
Poisson point processes. There one wanted the $\pi$ to be a sequence of 
{\it coarser and coarser} partitions, and also one needed some extra 
properties for the distribution of configuration points in the cells of 
$\pi$, which makes the proof lengthier (and restriction to Poisson point 
processes somehow necessary).

Fix a point $o$ of the plane.

\procl l.partition
Let $\omega$ be a point process such that the only isometry for the 
configuration is the identity a.s.. Let $\epsilon_i \to 0$ be 
arbitrary. Then there is a sequence of partitions 
$\pi$ of the plane, defined as equivariant functions of $\omega$, and such 
that 
$$\P[o\in C, C\in\pi, C \hbox{ is a } 2^i\times 2^i \hbox{ square}]\geq 
1-\epsilon_i.$$
\endprocl

\proof
Choose an equivariant subset $\omega_n\subset \omega$ such that any two 
elements of $\omega_n$ are at distance $\geq n $ from each other. See 
Corollary 3.2 in \ref 
b.T/ for such a choice. Let ${\cal V}_n$ be the Voronoi tessellation 
on $\omega_n$. Then, as shown in \ref b.T/, the probability that 
a point $x$ is in the $r$-neighborhood of the boundary of some cell in 
${\cal V}_n$ is at most $cr/n$ with some universal constant $c$. 
Now, let $n(i)$ be a sequence of integers that tends to infinity fast 
enough, and for each $C\in {\cal V}_{n(i)}$, subdivide $C$ by a copy of 
the 
$2^i\times 2^i$ square grid, whose position is determined by some 
deterministic 
rule 
(which tells, e.g., in which corner of $C$ one should put the origin of
the 
grid, and which incident edge should be ``covered" by the horizontal axis 
of the grid). Let the set of cells resulting from this subdivision be 
$\pi$. We have that $\P[o\in C, C\in \pi, C \hbox{ is not a } 2^i\times 
2^i \hbox{ square}]\leq \P[o$ is in the $2^{i+1}$-neighborhood of the 
boundary of some cell in
${\cal V}_{n(i)}]\leq c2^i/n(i)$. This is arbitrarily small, if $n(i)$ 
grows 
fast enough, proving the claim.  \Qed


The next example shows a translation invariant random planar map that does 
not have the regular decay 
property:
\procl x.aberralt
For simplicity, we construct a partition of $\Z ^d$ that is invariant with 
respect to translations of $Z^d$. One can easily modify this by random 
rotations to get an isometry-invariant partition of the plane.

For each $i\in \Z$, let $\xi_i$ be a geometric random variable with 
parameter $1/2$. Partition the vertical line $\{(i,j)\, :\, j\in\Z\}$ to 
intervals of length $2^{2^{\xi_i}}$ each, by choosing one of the 
$2^{2^{\xi_i}}$ 
such partitions uniformly, independently for the different $i$'s.

Similar but more complicated constructions lead to examples 
that are invariant 
under planar isometries, and look ``more 2-dimensional".
\endprocl

\procl p.voronoi
Let $\omega$ be a point process. 
Then the graph $G$ defined on $\omega$ by the Voronoi tessellation 
satisfies the 
regular decay property.

In particular, the Poisson-Voronoi map has the regular 
decay 
property.\endprocl

\proofof p.voronoi
Suppose that the statement is false. Then there is an $a>0$ such that for 
every $r$ there is an $R\geq r$ such that with probability at least $a$, 
$B_R(o)$ contains
$x,y\in \omega$ with adjacent Voronoi cells and such that $\dist_{\R^2} 
(x,y)\geq aR/6$.
Now, consider the square $S$ over diagonal $xy$ and the two triangles that 
the 
diagonal $xy$ divides $S$ into. It is easy to check that if both 
these triangles contain a configuration point in their interiors, then $x$ 
and $y$ cannot have adjacent Voronoi cells. So one of them has to be 
empty, consequently $S$ contains an empty square of diagonal half of 
that of $S$. 
We conclude that the probability that 
$B_R (o)$ 
contains $x,y\in \omega$ with adjacent Voronoi cells and such that 
$\dist_{\R^2}
(x,y)\geq aR/6$, is smaller than the probability that it 
contains an empty square $D$ of area $(Ra/6)^2 /4$. Covering  
$B_R$ 
of $o$ by $c a^2$ many squares of area $(Ra/6)^2/16$, one of them thus has 
to 
be empty (one that is inside $D$). Summing up 
the probabilities for this, we get
$\P[B_R$ contains a pair of adjecent vertices at distance $aR/6]\leq 
c a^{-2} 
\P[$a fixed square of area $(aR)^2 /576$ is empty$]$.  
Note that $c$ was a constant independent of $r$ and $R$, so this latter 
tends to 0 as $R$ tends to infinity. This contradicts the assumption on 
$a$.
\Qed

\procl p.fontos
If $G$ has the regular decay property, then there is a cycle exhaustion of 
width 6 for $G$.\endprocl

\def\pi{{\cal P}_i}
\proof
As before, $o$ is a point of the plane.

Let $\pi$ be a sequence of partitions of the plane such that $\P [o\in 
C,\, C\in\pi$ is 
an $r_i\times r_i$ square$]\geq 1-2^i$, as given by \ref p.voronoi/ 
setting $\epsilon_i=2^{-i}$ for simplicity. The $r_i$ will be chosen 
later, to increase fast enough.
Let $E_i$ be the set of edges in 
$G$ 
that 
intersect the boundary of some cell in $\pi$, and let $E_i^j$ be the set 
of edges in $G$ at distance $\leq j$ from $E_i$ (hence $E_i^0$ is $E_i$, 
$E_i^1$ is the set of 
edges of $G$ with an endpoint in $E_i$, etc). 
Let $G_i :=G\setminus E_i^4$.

For a subset $A$ of the plane, let $\partial_r A $ be the set of point at 
Euclidean distance at most $r$ from $A$.
Say that $C\in\pi$ is {\bf good}, if it is an $r_i\times r_i$ square and 
there is no path of length $\leq 4$ in $G$ that connects the complement of 
$C$ 
with $C^o :=C\setminus \partial_{a(r_i) r_i} C$. Here $a(r)$ is the 
function from the definition of the regular decay property.
By the assumption on 
$\pi$ 
and using the definition of the positive decay property, we obtain 
$$\P[x\in C, C\in\pi \hbox{ is good} ]\geq 1-2^{-i}-a(r_i). \label e.1
$$

Now, if $C$ is good, then all vertices in $C^o$ are contained in the same
connected component of $G_i$: otherwise the graph induced by $E_i^4\cup 
(G\setminus C)$ would separate them, which implies that some edge of 
$E_i^5$ 
would cross $C^o$. Then there would be a path of length at most 11 
containing 
this edge and crossing the boundary of $C$ by both its first and last 
edge; in particular one of the edges in this path would have length $\geq 
2a(r_i) r_i/11> a(r_i) r_i/6$, contradicting the assumption that $C$ is 
good. We have 
obtained 
that for a fixed point x of the plane:
$$\P[x\in C^o, C\in \pi, C^o\cap V(G) \;\hbox{is in one connected}$$ 
$$\hbox{component of } G_i]\geq 1-2^i-3a(r_i) \label e.fontos
$$
using \ref e.1/and the generous upper bound $2a(r_i)$ on the probability 
that $x\in C\setminus C^o$.
From this it is easy to see that the probability that $x$ is surrounded 
by a cycle of $G_i$ also tends to 1 as $i$ tends to infinity.

\def\gi{G_i^{\rm good}}
Note that by definition every component $\gamma$ of $G_i$ is inside some 
set (cell) of the partition $\pi$. Call this 
$C(\gamma)$. Take $\gi$ 
to be the union of connected components $\gamma$ of $G_i$ such 
that 
every vertex inside $\gamma\cap C(\gamma)^o$ is in the same 
component of $G_i$. By \ref e.fontos/ and the remark after it we know 
that $\P[x$ is surrounded by a cycle in $\gi]$ tends to 1 with $i$.
By definition of $G_i$, every two connected components of $\gi$ have 
distance at least 8 (the 4-neighborhood of $E_i$ is in between two such 
components). Hence, for $i$ fixed, the set $B_i$ of external boundaries  
of the components of $\gi$ as in \ref p.boundary/ forms a family of 
non-selftouching cycles at distances at least 6 from each other. Observe 
that 
every cycle in $B_i$ is contained in $C\setminus C^o$ for some good $C\in 
\pi$, since it is in the boundary of a graph that contains $C^o\cap G$, 
but does not contain any element of $E_i^4$. We have seen that $x$ is 
surrounded by one cycle of $B_i$ with probability arbitrary close to 1 if 
$r_i$ is large enough. Another consequence of that every cycle 
$O\in B_i$ is 
contained in some $C(O)\setminus C(O)^o$, $C(O)\in\pi$ is the next 
assertion. For $j>i$, $O_j\in B_j, O_i\in 
B_i$, $C(O_i)$ can intersect the 5-neighborhood of $C( O_j)\cap G$ in $G$ 
only 
if the Euclidean distance of
$C(O_i)$ from the boundary of $C( O_j)$ is less then $5a(r_j) r_j$ (using 
that $C(O_j)$ is good). If 
$r_j$ was chosen to grow fast enough, the probability that the $C(O_i)$ 
containing $x$ is such for some $j>i$ tends to zero. That is, if we delete 
every $O$ with this property, then the probability that $o$ is 
contained in some cycle $O\in B_i$ that was not deleted, tends to 1 with 
$i$ arbitrarily fast by a suitable choice of $r_i$. Hence we can finish 
the construction as described in the next 
paragraph.

Delete every cycle of $B_i$ that intersects the 
5-neighborhood (in $G$) of any cycle in ${\cal P}_j$, $j>i$ arbitrary. 
Call the set of remaining cycles $\tilde B_i$. If the $r_i$ grew fast 
enough, the probability that the cycle of $B_i$ surrounding $x$ 
(conditioned on that there is such a cycle) intersects the cycle of 
some ${\cal P}_j$, $j>i$, is at most $2^{-i}$. 
The 
probability that a cycle of $\tilde B_i$ surrounds $x$ is at least 
$1-2^{-i}$.

We conclude that $\cup \tilde B_i$ is a cycle exhaustion. The 
corresponding tree $T$ and labelling of the vertices of $T$, is uniquely 
determined by the 
construction (see the comment after the definition of a cycle exhaustion).
\Qed

\procl l.exhaustion
Let $G$ be a random triangulated planar map that satisfies the
regular decay property, and suppose
that $G$ has only the trivial symmetry a.s.. Then there
exists an equivariant function of $G$ that is an even
cycle exhaustion of corridor width $4$.
\endprocl

\proof
To prove the existence of a cycle exhaustion with even cycles can be 
obtained as a modification of the cycle exhaustion constructed in \ref 
p.fontos/. Note that if the   
set $\{v\in V(T)\,: \,\lambda (v)$ is even$\}$ 
has
a complement in $T$ with only finite components, then keeping only the 
even
$\lambda (v)$'s, we would obtain an even cycle exhaustion. Hence, if this 
is not the case, one may keep only the odd cycles and get a cycle 
exhaustion. So, consider a cycle exhaustion with only odd cycles.
Call the set of cycles corresponding to the leaves of the tree in 
the cycle exhaustion ${\cal L}_1$, those corresponding to neighbors of 
the leaves that are not leaves ${\cal L}_2$, and so on. Call the set of 
cycles in 
the exhaustion ${\cal C}$. That is, ${\cal C}=\{\lambda (x)\, :\, x\in 
V(T)\}$. We will keep notation $\int (O)$ when $O\in {\cal C}$, to denote 
the bounded component of $G\setminus O$.
Now, we will show that one is able to modify any cycle $O_1 \in {\cal C}$ 
and some $O_0\in {\cal C}$ inside $\int (O_1)$, 
to 
get 
an even cycle $\nu (O_1)$, preserving the 
property 
that the $\nu (O_1)$'s are at distance 
at least $4$ from each other. 

Let $O_1$ be an arbitrary odd cycle in ${\cal C}$, such that there is an 
$O_0\in {\cal C}$ 
contained in $\int (O_1)$, chosen in a later 
defined way. We will find a way to remove a small arch of $O_1$, and 
connect the remaining arch of $O_1$ to an arch of $O_0$ by two paths in 
such a way that the resulting graph is a non-selftouching even cycle, it 
still has 
distance $\geq 4$ from the other cycles of ${\cal C}$ (or their modified 
version, if we have already modified them in the way we are modifying 
$O_1$), and further, it surrounds ``almost" as many points as $O_1$ did, 
so 
condition 
(5) of a cycle exhaustion is preserved by the 
modified 
cycles. (See Figure 4. for an illustration of what follows.) More 
precisely, we will find:
\item{(I)} Paths $P_1$ and $P_2$ in $\int (O_1) \setminus O_0$, such that 
there is an endpoint $x_i$ for $P_i$ that is adjacent to $O_0$, the other 
endpoint $y_i$ of $P_i$ is adjacent to $O_1$, and the number of vertices 
in $O_0$ that are adjacent to $x_0$ and to $x_1$, respectively, have the 
same parity.

\item{(II)} $P_1$ and $P_2$ are not self-touching, they do not touch each 
other, and none of their inner vertices 
is 
adjacent to $O_0\cup O_1$.

\item{(III)} $P_1\cup P_2$ has distance $\geq 4$ from all $\bar 
O\in{\cal 
C}\setminus \{O_0, O_1\}$.

\item{(IV)} Every child $O'\not = O_0$ of $O_1$ is in the same connected 
component of $G\setminus (P_1\cup P_2\cup O_1\cup O_2)$.

Suppose we can find the above described objects. Let $\ell_1$ and $r_1$ 
($\ell_2$ and $r_2$) be the ``extremal" neighbors 
of $P_1$ ($P_2$) on $O_0$. By extremal we mean that there is no neighbor 
of $x_1$ in one of the archs of $O_0$ from $\ell_1$ to $r_1$ (and 
similarly for $x_2$). Index them so that 
the cyclic order of these 4 points on 
the cycle $O_0$ is $\ell_1,r_1,\ell_2,r_2$. Let the arch between $r_1$ and 
$\ell_2$  (respectively $r_2$ and $\ell_1$) that does not contain the 
other two points be $A_1$ (respectively $A_2$). Finally, let $A$ be the 
longer of the two archs on $O_1$ between a neighbor of $P_1$ and a 
neighbor of $P_2$ such
that $A$ does not contain any other neighbor of $P_1\cup P_2$.
Note that $(-1)^{|A_1|+|A_2|}=(-1)^{|O_0|}=-1$, where the first equation 
is by  
(I) and the second is by the assumption that every cycle in ${\cal C}$ is 
odd. Hence one of $A\cup 
P_1\cup 
P_2\cup A_1$ and $A\cup P_1\cup  
P_2\cup A_2$ is even (since they have opposite parity); call this $\nu 
(O_1)$. 
Note that $\nu (O_1)$ is a 
non-selftouching cycle, by (II), (III) and the assumption that $O_0$ and 
$O_1$ had distance at least 6. 
Now, if we consider $\nu (O_1)$ for every cycle
$O_1\in{\cal L}_{2k}$, $k\in\Z^+$, then no cycle of ${\cal C}$ is used as 
$O_1$ or $O_0$
for more than one $\nu (O_1)$. On the other 
hand, condition (IV) 
guaranteed that the 
interior of $\nu
(O_1)$ contains every element of ${\cal C}\setminus \{O_0\}$ that $\int
(O_1)$ contained, so (5) remains valid for $\{\nu (O_1)\}$.
Hence the resulting set $\{\nu (O_1)\, :\, O_1\in{\cal L}_{2k}, 
\,k\in\Z^+\}$
is a cycle exhaustion. 

Therefore it only remains to show the existence of $x_1,x_2$ and $P_1,P_2$ 
that 
satisfy (I)-(IV). 
So let $O_1\in {\cal C}$ be given, and let $Q$ be a non-selftouching path 
with endpoints adjacent to $O_1$ and $O_0$ respectively,
where $O_0\subset \int (O_1)$ is chosen so that $Q\subset 
G\setminus\cup_{O\in{\cal C}, O\not= O_0}(O\bigcup \partial_{c-2} 
O)$. By switching to a subpath of $Q$ if necessary, we may assume that 
none of the inner vertices of $Q$ is adjacent to $O_1\cup O_0$.
If we 
take the outer boundary {\it of} $\partial Q$ in $G$, then it contains two 
non-selftouching paths $Q_1$ and $Q_2$, between $O_0$ and $O_1$. (Make 
them non-selftouching by choosing them to have minimal length.)
By switching to subpaths of $Q_1$ and $Q_0$ if necessary, we may assume 
that    
none of the inner vertices of $Q_1$ or $Q_0$ is adjacent to $O_1\cup O_0$.
Since the parity condition in (I) has 
to be satisfied by at least two of $Q_1,Q_2,Q$, we can chose those two 
to be $P_1$ and $P_2$. One can easily check that the other requirements 
are also satisfied. 
\Qed
$$$$
\SetLabels
  (.5*-.17) {{\bf Figure 4.} The construction of $\nu (O_1)$. Here 
$P_1:=Q_1$ and $P_2:=Q_2$.}\\
  (.23*.89) {$Q_2$}\\
  (.26*.90) {$Q$}\\
  (.33*.87) {$Q_1$}\\
  (.22*.77) {$x_2$}\\
  (.33*.76) {$x_1$}\\
  (.07*.9) {$O_1$}\\
  (.35*.57) {$O_0$}\\
  (.95*.9) {$\nu (O_1)$}\\
\endSetLabels
\AffixLabels{\epsfysize=4.5cm
\epsfbox{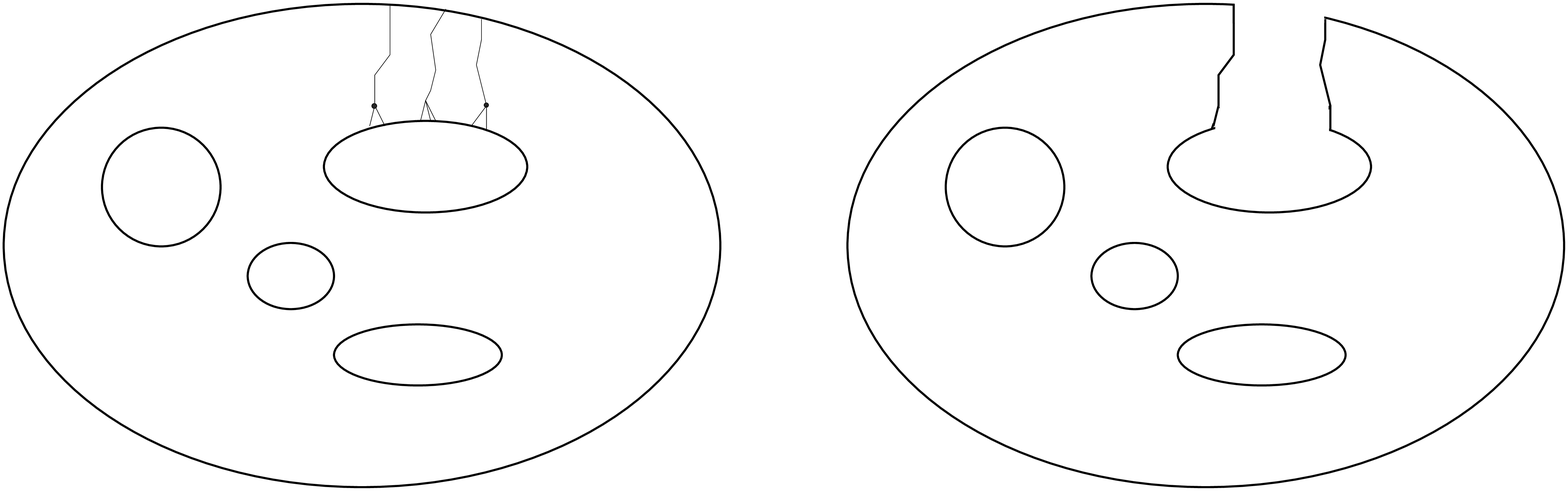}}
$$$$
\bigskip

\bsection{Concluding remarks, further directions}{s.4}

In this section we discuss how necessary the conditions in \ref t.main/ 
are. We characterize the case when there is some nontrivial symmetry. 
Whether the conclusion of \ref t.main/ holds when we do not assume the 
regular decay property is not clear. By \ref l.exhaustion/ it would follow 
from a positive answer to the next question. 

\procl q.cycle 
Let $G$ be a random graph in the plane, whose distribution is invariant
with respect to some transitive group of isometries of the plane. Suppose 
that almost surely the only symmetry that $G$ has is the trivial one. 
Is there an even cycle exhaustion of corridor width 4 for $G$?
\endprocl


\procl r.generalizations
We assumed that $G$ has only trivial symmetries. This assumption can be
slightly
weakened, since we only need the lack of symmetries in order to
construct the sequence of partitions ${\cal P}_i$ in the plane, and to
make ``local choices"...
So, suppose that $G$ is an arbitrary graph embedded in the
plane, with an invariant distribution, and that there is some invariant
point
process $P$, which may not be independent of $G$. Then one may look at 
colorings of $G$ that are equivariant measurable functions of {\it the 
pair} 
$(P,G)$.  
We usually assume 
that one of $P$ 
and
$G$ is an equivariant function of the other.
One type of example
is when $G$ is an equivariant function of $P$, such as the graph
given by the Voronoi map
on $P$; another class of examples is when we a priori have
$G$, and then define $P$ as an equivariant function of this, e.g. $P$
is the set of vertices in $G$.
If $P$ and $V$ are independent, that correspond to the case when
one can use local extra randomness while coloring $G$. See the next 
example as an illustration of this more general setup (and a case 
where there is no equivariant 4-coloring).
\endprocl

\procl x.uj
Let $\omega$ be the point process obtained as follows. Let $H$ be
the triangular grid of unit edge lengths, with an extra vertex added in 
the center of each triangle and connected to the three nodes of the 
triangle.
Translate $H$ by a uniformly 
chosen vector from the 
union of the 6 triangles of the original triangular lattice incident to 
some vertex. We get a point set 
$\omega '$ that 
is invariant with respect to translations. Now, relocate every point of 
$\omega$ uniformly in its neighborhood of radius $1/100$. The resulting 
$\omega$ has only trivial isometry almost surely, and the Voronoi 
tessellation on $\omega$ as a map is isomorphic to $H$. Hence its 
chromatic number is 4, and our method gives a 5-coloring that 
is an equivariant measurable
function of $\omega$. On the other hand, 4 colors do not suffice for this, 
since up to permutation of colors there is a unique 
4-coloring for $H$.
\endprocl 


\procl r.lattice
Consider now the general case when $G$ does have some nontrivial symmetry 
with 
positive probability. 

Then each ergodic component where $G$ has some nontrivial symmetry is the 
$\Gamma$-translate of some 
quasi-transitive graph $H$. Consider $H/\Gamma$. The subgroup of $\Gamma$ 
of elements whose natural action on the torus defines an automorphism for 
the $H/\Gamma$ embedded in the torus is trivial. Hence, if $H/\Gamma$ is 
colorable by $k$ 
colors, then that extends to an equivariant measurable coloring of $H$. 
We get the color of each $x\in G$ by simply identifying   
which vertex of $H/\Gamma$ the factor map maps $x$ into (which can be 
determined from a large enough neighborhood of $x$), and taking the 
color of that vertex by $c$.
Otherwise there is no 
coloring by any number of colors (this is the case when $H$ has a loop edge). 

Conversely, any finite graph $F$ embedded in the torus can be lifted to 
define 
a 
quasi-transitive graph $H$ embedded in the plane, and a random translate 
can be used to define 
$G$. The chromatic number of $F$, is either between $1$ and $7$ or $F$ is not 
colorable by any number of colors, see \ref 
b.ThC/. Hence, if $F$ is embedded in the torus so that there is a coloring by 
$k$ colors such that 
$H$ is 
also equivariantly $k$-colorable, otherwise it is not.
 
We have obtained the following:

\procl t.rest
Let $G$ be an ergodic random graph in the plane, whose distribution is 
invariant   
with respect to some transitive group of isometries of the plane, and 
suppose that $\Gamma$ acts quasitransitively on $G$. 
Define 
$\Delta:=G/\Gamma$ (which belong to one graph-isomorphism class almost 
surely). 
Then the minimal number of colors needed for an 
equivariant coloring of $G$ is the chromatic number of $\Delta$. This can 
be any number between 1 and 7; or, if $\Delta$ has a loop-edge, then $G$ 
is not colorable in an equivariant measurable way by any number of 
colors. \endprocl

\procl q.4
Is it true that for any random planar graph $G$ that is invariant with 
respect to some transitive group $\Gamma$ of isometries of the plane, and 
that has no nontrivial 
symmetries, there is a 4-coloring as an equivariant measurable function of 
$G$?

A necessary condition for a positive answer is that $G$ has infinitely 
many 4-colorings, which is, to our knowledge, also open in graph theory.

Let us call a partition of the vertex set of a graph $G$ into classes
$\{K_1,\ldots ,
K_k\}$ such that each
of the $K_i$ is an independent set a {\bf blind
k-coloring}. Note that the color classes of every $k$-coloring give rise
to a
blind $k$-coloring. However, there is an invariant graph in 1 dimension
(actually, every invariant connected graph, a biinfinite path, is such)
that has an invariant blind 2-coloring, but no invariant 2-coloring.
Take e.g. a Poisson point process on the line, and let $G$ be the
corresponding Voronoi map. There is no way to 2-color $G$ in an
equivariant
way, because that would contradict ergodicity of the point process.
On
the other hand, there is a trivial equivariant blind 2-coloring: let the
interval of     
the origin and all intervals at an even distance from it form $K_1$, and
the other intervals $K_2$. The way we chose $K_1$ is of course not
invariant, but the {\it set} $\{K_1,K_2\}$ is. Perhaps surprisingly, it
would be
a
lot easier to show the existence of a blind 5-coloring for the case of
\ref t.main/,
than it was to show the existence of a 5-coloring, for the reason the we
sketch in the next few paragraphs.

The last part of the proof of \ref l.exhaustion/ consisted of showing that
one can find a
cyclic
exhaustion {\it
consisting of even cycles}. If we were satisfied with a blind 5-coloring,
this last step could be omitted by some modification of \ref l.cycles/:

\procl l.blind
Let $(T,\lambda)$ be a cycle exhaustion of $G$. Then there is an
equivariant blind 5-coloring of $G$.
\endprocl

The proof proceeds similarly to that of \ref l.cycles/, with the following
differences. When we define the $G_v$ we contract all {\it but one}
vertex
of each odd $\lambda (v)$. Then, for $w\to v$, we ``match" the
colorings
of $\lambda (w)$ determined by $\gamma_v$ and by $\gamma_w$ (in a way
defined shortly) by permuting
the colors assigned by $\gamma_w$. This means infinitely many permutations 
of
colors as
$v\in V(T)$ goes to infinity, hence we loose the color and can only
detect whether two points are in the same or different color classes.
This is exactly a blind 5-coloring. The only thing missing from this 
sketch is how $\gamma_v$ would tell the color of an odd $\lambda
(w)$
(before any potential permutations): color every second vertex, and the
vertex that was not contracted, similarly to what $\gamma_v$ colored   
them (or their image after the identification), and color the remaining
vertices with 0.

There are some questions of similar flavor that we would like to mention 
to finish with. The first one is purely deterministic.

\procl q.BL
Does every quasi-transitive planar graph $G$ admit a periodic
4-coloring?

This was first asked by Bowen and Lyons. They observed that when the graph 
is Euclidean, there exists a 5-coloring, by the following argument. It was 
shown by Thomassen \ref b.ThC/ that every graph embedded on a surface of 
genus $g>0$ with all noncontractible cycles long enough, can be colored by 
5 
colors. 
If $\Gamma$, as before, is the
fixed transitive
group of
isometries of the plane that also act as
automorphisms of $G$, the quotient of the plane by $\Gamma$ is a torus, 
with an embedded graph $H=G/\Gamma$.
Any coloring of $H$ lifts to a periodic coloring 
of $G$. We may assume that the noncontractible cycles of $H$ are long 
enough for the assumption of Thomassen's theorem (and hence the 
conclusion that $G$ has a periodic 5-coloring), otherwise replace $H$ 
by the $H'$ obtained when lifting $H$ to a torus that covers that of $H$ 
$k$ times ($k$ large enough). 

\procl q.LS
Does every infinite quasi-transitive graph have an invariant random
coloring with as many colors as its chromatic number? 

A graph is called quasi-transitive, if its set of automorphisms have 
finitely many orbits on the vertices. Invariance of a coloring is 
understood with respect to this group.
This question was asked by Lyons and 
Schramm. The latter has shown that a coloring by $d+1$ colors, where $d$ 
is 
the maximal degree in the graphs, is always possible.

\procl q.L
How many colors are needed to have a mixing invariant random
coloring? What if the coloring has to be an equivariant function of i.i.d.
${\rm Unif}[0,1]$ random labels on the vertices?

This question is originated from Lyons. It is easy to see (by arguments 
similar to what we used when talking about blind-colorability and 
colorability in one dimension) that transitive trees can be invariantly 
colored by 2 colors, but one needs 3 to get a mixing 2-coloring.

\medbreak

\startbib                                
  
\bibitem[1]{ABGMP} Angel, Benjamini, Gurel-Gurevich, Mayerovitch \and 
Peled (2008) Stationary map coloring ({\it preprint}).

\bibitem[2]{Th} Thomas, R. (1998) An update on the Four-color Theorem
{\it Notices Amer. Math. Soc.} {\bf 45}, no. 7, 848--859.

\bibitem[3]{ThC} Thomassen, C. (1993) Five-coloring maps on surfaces 
{\it J. Combin. Theory Ser. B} {\bf 59}, no. 1, 89--105.

\bibitem[4]{T}
Tim\'ar, \'A. (2004)  Tree and grid factors for general point
processes.
{\it Elec. Comm.
in Probab.} {\bf 9}, 53--59.

\bibitem[5]{T2} Tim\'ar, \'A. (2008) Invariant matchings of exponential
tail on coin flips in $\Z^d$ ({\it preprint}).

\endbib 
\bibfile{\jobname}
\def\noop#1{\relax}
\input \jobname.bbl

\filbreak
\begingroup
\eightpoint\sc
\parindent=0pt\baselineskip=10pt

Hausdorff Center for Mathematics, Universit\"at Bonn, D-53115 Bonn

\emailwww{adam.timar[at]hcm.uni-bonn.de}{}
\htmlref{}{http://www.hausdorff-center.uni-bonn.de/people/timar/}
\endgroup

\bye